\documentclass[12pt,a4paper,twoside]{article}
\usepackage[latin1]{inputenc}
\usepackage{amsmath,amsfonts,amssymb,a4,enumerate,epsfig,array,theorem}

\def\qed{\unskip\nobreak\hfil\penalty50\hskip1.75em\null\nobreak\hfil
$\blacksquare$ {\parfillskip=0pt \finalhyphendemerits=0 \par}\goodbreak}

\newfont{\eightrm}{cmr8}
\newfont{\fiverm}{cmr5}

\newcommand\nobf{\noindent\bf}
\newcommand\proof{{\noindent\bf Proof: }}

\newcommand{\ZZ}{{\mathbb{Z}}}

\newcommand{\RR}{{\mathbb{R}}}

\newcommand{\Ss}{{\mathbb{S}}}
\newcommand{\Ff}{{\cal F}}
\newcommand{\Jj}{{\cal J}}
\newcommand{\rot}{{\tau}}
\newcommand{\id}{\textrm{id}}
\newcommand{\Diff}{\textit{Diff}\,}

\newtheorem{theo}{Theorem}
\newtheorem{prop}{Proposition}[section]

\newtheorem{lemma}[prop]{Lemma}

\begin{document}
\title{Nilpotent pseudogroups of \\ functions on an interval}
\author{Tania M. Begazo and Nicolau C. Saldanha}
\maketitle

\begin{abstract}
A {\it near-identity} nilpotent pseudogroup of order $m \ge 1$ is
a family $f_1, \ldots, f_n: (-1,1) \to \RR$ of $C^2$ functions for which:
$|f_i - \id|_{C^1} < \epsilon$ for some
small positive real number $\epsilon < 1/10^{m+1}$ and
commutators of the functions $f_i$ of order at least $m$
equal the identity.
We present a classification of near-identity nilpotent pseudogroups:
our results are similar to those of Plante, Thurston, Farb and Franks.
As an application, we classify certain foliations of nilpotent manifolds.
\end{abstract}


{\noindent\bf MSC-class:} 37E05; 37E45; 57S25; 57R30

\section{Introduction}

A {\it near-identity} nilpotent pseudogroup of order $m \ge 1$ is
a family $f_1, \ldots, f_n: (-1,1) \to \RR$ of $C^2$ functions for which:
\begin{itemize}
\item{$|f_i - \id|_{C^1} < \epsilon$ for some
small positive real number $\epsilon < 1/10^{m+1}$;}
\item{commutators of the functions $f_i$ of
order at least $m$ equal the identity.}
\end{itemize}
More precisely, the pseudogroup is generated by the family of functions
but we shall often make the abuse of confusing these two objects.
A commutator of order $1$ is of the form $[f_i,f_j]$
and, for $m > 1$, a commutator of order $m$ is a function of the form
$[g_1,g_2] = g_1^{-1} \circ g_2^{-1} \circ g_1 \circ g_2$
where $g_1$ is one of the original functions or
a commutator of order less than $m$ and
$g_2$ is a commutator of order $m-1$.
In particular, a family of functions which commute is
a nilpotent pseudogroup of order 1: we call this an {\it abelian} pseudogroup.
If all commutators commute with each other,
then we call the pseudogroup {\it metabelian}.
A common fixed point for the pseudogroup is a point $x$ for which $f_i(x) = x$
for all $i$.

If the functions $f_i$ are bijections from $(-1,1)$ to itself, then
such a nilpotent pseudogroup is just a group of functions.
Plante and Thurston prove that nilpotent groups of diffeomorphisms
on $[-1,1)$ or $(-1,1]$ are abelian (\cite{PT}).
Farb and Franks consider groups of diffeomorphisms of $(-1,1)$
and prove several results, among them that $\Diff^\infty((-1,1))$
contains nilpotent subgroups of arbitrary
degree of nilpotency but that such subgroups are all metabelian (\cite{FF}).
Our work differs from these in that we consider pseudogroups
instead of groups: the main new difficulty is that long compositions
of the functions $f_i$ may not define any function
since the domain may vanish.
More precisely, we prove the following theorem.

\begin{theo}
\label{theo:classify}
Any near-identity nilpotent pseudogroup
of functions $f_1, \ldots, f_n$ is metabelian.
Furthermore, any near-identity pseudogroup
fits into one of the three cases below.
\begin{enumerate}
\item{There exists at least one common fixed point and
the pseudogroup is abelian. Furthermore, for each maximal
interval $I_1 \subset (-1,1)$ containing no common fixed points,
there exist real constants $a_i$
and an increasing homeomorphism $\phi: J \to I_1$, $J \subseteq \RR$ with 
$f_i(\phi(t)) = \phi(t + a_i)$ whenever $t, t+a_i \in J$.
If $\inf I_1 > -1$ (resp. $\sup I_1 < 1$)
then $\inf J = -\infty$ (resp. $\sup J = +\infty$).}
\item{There exists no common fixed point, the pseudogroup is abelian,
there exist real constants $a_i$
and an increasing homeomorphism $\phi: J \to (-1,1)$, $J \subseteq \RR$ with 
$|J| > |a_i|$ and $f_i(\phi(t)) = \phi(t + a_i)$ whenever $t, t+a_i \in J$.}
\item{There exists no common fixed point,
there exist integer constants $a_i$ and
a finite set $\{y_{-N}, y_{-N+1}, \ldots, y_N\} \subset (-1,1)$
with $y_i < y_{i+1}$, $y_{-N} < -1+\epsilon$, $y_N > 1-\epsilon$,
$N > |a_i|$, such that $f_i(y_k) = y_{k + a_i}$.}
\end{enumerate}
\end{theo}

This subject was motivated by the study of actions of nilpotent Lie groups
on manifolds. In particular, we study actions of the Heisenberg group
\[ G = \left\{
\begin{pmatrix} 1 & 0 & 0 \\ x & 1 & 0 \\ z & y & 1 \end{pmatrix},
x,y,z \in \RR \right\} \]
on manifolds of dimension 4 (\cite{B}, \cite{BS}).
We do not discuss actions here but we present an application
of theorem \ref{theo:classify} to foliations.

We consider compact orientable manifolds of the form $G/H = \{ gH, g \in G\}$,
where $G$ is a nilpotent Lie group
and $H = \pi_1(G/H) \subset G$ is a discrete cocompact subgroup.
We may assume that $G/H$ has a smooth Riemann metric
and therefore $G/H \times (-1,1)$ also has a metric.
The foliation $\Ff_0$ of $G/H \times (-1,1)$ with leaves of the form
$G/H \times \{x\}$ is called {\it horizontal};
this foliation can be defined as being perpendicular
to the vertical vector field $Z_0$ at any point of $G/H \times (-1,1)$.
An arbitrary transversally orientable foliation $\Ff_1$
of $G/H \times (-1,1)$ of codimension $1$ can be similarly described
as being perpendicular to some unit vector field $Z_1$ at any point;
we define $d_{C^1}(\Ff_0,\Ff_1) = d_{C^1}(Z_0,Z_1)$.
Let $\gamma: [0,1] \to G/H$ be a smooth path;
if $d_{C^1}(\Ff_0,\Ff_1)$ is sufficiently small
then, for $x \in (-1+\epsilon, 1-\epsilon)$,
the path $\gamma$ can be {\it lifted} to
a unique smooth path $\gamma_x: [0,1] \to G/H \times (-1,1)$
with $\gamma_x(0) = (\gamma(0),x)$, $\gamma_x(t) = (\gamma(t), \ast)$
and $\gamma_x'(t)$ tangent to $\Ff_1$ for all $t$  (\cite{Camacho}).
This defines a function $f_\gamma: (-1+\epsilon, 1-\epsilon) \to (-1,1)$
taking $x$ to the second coordinate of $\gamma_x(1)$:
if $\gamma$ is a generator of $\pi_1(G/H)$,
$f_\gamma$ is called the {\it holonomy} of $\Ff_1$.
If $\gamma_1, \gamma_2, \ldots, \gamma_n$ are generators of $\pi_1(G/H)$
and $d_{C^1}(\Ff_0, \Ff_1)$ is sufficiently small 
then $f_{\gamma_1}, f_{\gamma_2}, \ldots, f_{\gamma_n}$ form
a near-identity nilpotent pseudogroup.
We call a foliation of $G/H \times (-1,1)$ {\it abelian}
(resp. {\it metabelian}) if its pseudogroup is abelian (resp. metabelian).
Similarly, we call a leaf of a foliation {\it abelian}
if $\gamma_x(1) = \gamma_x(0)$ whenever $\gamma \in \pi_1(G/H)$
is a commutator and $\gamma_x(0)$ belongs to the leaf;
thus, a foliation is abelian if and only if all its leaves are abelian.

Let $\alpha: G \to \RR$ be a group homomorphism:
the $1$-form $d\alpha$ can be lifted to $G/H$
and therefore to $G/H \times \RR$,
where it is closed but probably not exact.
The $1$-form $dx$, where $x$ is the second coordinate,
is exact in $G/H \times \RR$ and therefore $d\alpha + dx$ is a closed
$1$-form in $G/H \times \RR$, defining a foliation $\Ff_\alpha$.
Notice that $\Ff_\alpha$ is abelian.

\begin{theo}
\label{theo:foliation}
Given a manifold of the form $G/H$ there exists $\epsilon_{G/H}$ such that,
if $\Ff_1$ is a smooth foliation of $G/H \times (-1,1)$
with $d_{C^1}(\Ff_0,\Ff_1) < \epsilon_{G/H}$,
then $\Ff_1$ is metabelian and one of the three conditions hold:
\begin{enumerate}
\item{$\Ff_1$ has at least one compact leaf, is abelian and,
for any maximal connected open set $M$
of $G/H \times (-1,1)$ containing no compact leaf,
there exists a homomorphism $\alpha: G \to \RR$,
an open set $\Jj \subseteq G/H \times \RR$
and a homeomorphism $\Phi: \Jj \to M$
taking $\Ff_\alpha$ to $\Ff_1$.}
\item{$\Ff_1$ has no compact leaf, is abelian and
there exists a homomorphism $\alpha: G \to \RR$,
an open set $\Jj \subseteq G/H \times \RR$
and a homeomorphism $\Phi: \Jj \to G/H \times (-1,1)$
taking $\Ff_\alpha$ to $\Ff_1$.}
\item{$\Ff_1$ has no compact leaf
and has an abelian leaf closed in $G/H \times (-1,1)$,
arriving or accumulating both at $G/H \times \{1\}$ and $G/H \times \{-1\}$.}
\end{enumerate}
\end{theo}

Some key results, presented in section 2, are that there exist
nonempty closed sets $X$ invariant under the pseudogroup of functions
and such that the restrictions to $X$
of all commutators equal the identity;
in particular, $[f_i,f_j]$ has many fixed points.
In section 3 we define a concept of {\it translation number} $\rot(f_i,f_j)$
for functions $f_i, f_j$ in a near-identity nilpotent pseudogroup;
this concept is also present in \cite{FF} and
reduces to the usual definition of rotation number if $[f_i,f_j] = \id$.
As we shall see in section 4,
Denjoy's theorem implies that if some translation number is irrational
then the invariant set $X$ constructed in section 2 is an interval.
In section 5 we use Koppel's lemma to
prove that near-identity nilpotent pseudogroups are metabelian
and finally, in section 6, we bring together the results
in order to prove theorems \ref{theo:classify} and \ref{theo:foliation}.

\medskip

We thank Carlos Gustavo Moreira for helpful conversations
and the referee for several helpful suggestions.
The authors acknowledge support from CNPq, CAPES and Faperj (Brazil).
This work was done while the first author was in a leave of absence
from Universidade Federal do Par\'a and had a visiting position in UFF.

\section{Fixed points}

We recall the usual definition of translation number for continuous
increasing functions $u: \RR \to \RR$ of degree one,
i.e., with $u(x+1) = u(x) + 1$, or, more generally,
for a continuous increasing function $u: [0,1] \to \RR$ with $u(1) = 1 + u(0)$.
Assume $0 < u(0) < 1$. Define the sequence $a_0 = 0$,
\[ a_{n+1} = \begin{cases} u(a_n),&\textrm{if } u(a_n) < 1,\\
u(a_n) - 1,&\textrm{otherwise.} \end{cases} \]
The translation number $\rot(u)$ of $u$ is the proportion of points
of this sequence in the interval $[0, u(0))$.
More precisely, define $p(0) = 0$, $p(n+1) = p(n) + 1$
if $0 \le a_n < u(0)$ and $p(n+1) = p(n)$ otherwise:
then the limit
\[ \lim_{n \to \infty}{\frac{p(n)}{n}} \]
exists and is called $\rot(u)$.
Denjoy's theorem \cite{D} states that if $u$ is a degree one function
of class $C^2$ with irrational translation number $\alpha = \rot(u)$
then there exists a homeomorphism $\phi: \RR \to \RR$, also of degree one,
such that $u(\phi(t)) = \phi(t + \alpha)$ for all $t$.

\medskip

The condition $|f_i - \id|_{C^1} < \epsilon < 1/10^{m+1}$
guarantees that $f_i$ is a diffeomorphism from $(-1,1)$ to
some open interval $I$,
\[ (-0.99, 0.99) \subset (-1+\epsilon, 1-\epsilon)
\subset I \subset (-1-\epsilon,1+\epsilon) \subset (-1.01, 1.01) \]
and we take $f_i^{-1}: I \cap (-1,1) \to \RR$.
A composition such as the commutator
\[ [f_i,f_j] = f_i^{-1} \circ f_j^{-1} \circ f_i \circ f_j \]
is defined in the largest possible domain
such that intermediate expressions are in $(-1,1)$:
this is clearly some interval $I$, $(-0.96, 0.96) \subset I \subseteq (-1,1)$.
An equality such as $[f_i,[f_j,f_k]] = \id$ means that both functions
coincide in the intersection of their domains,
in this case at least $(-1+10\epsilon,1-10\epsilon)$.
A set $X \subseteq (-1,1)$ is said to be invariant under the functions $f_i$
if $x \in X$ implies that $f_i(x)$ and $f_i^{-1}(x)$ belong
to $X \cup (-\infty,-1] \cup [1,\infty)$ for all $i$.
Recall that $X \subseteq (-1,1)$ is closed in $(-1,1)$
if $\overline{X} \cap (-1,1) = X$.

\begin{prop}
\label{prop:fixedpoint}
For any near-identity nilpotent pseudogroup of functions there exists
a nonempty set $X \subseteq (-1,1)$ which is closed in $(-1,1)$,
invariant under the functions $f_i$ and such that
$[f_i,f_j]|_X = \id$ for all $i,j$.
\end{prop}

If $x$ is a common fixed point, $X = \{x\}$ is a trivial
example of a closed invariant set as described in the proposition.
We shall be more interested in nontrivial invariant sets $X$.
Notice that if $x_0$ and $x_1$ are common fixed points
then the proposition may be applied to the restriction
of the functions to the interval $(x_0, x_1)$:
in other words, inside each maximal interval $I$ in the complement
of the set of common fixed points, there exists a nonempty invariant set
$X_I \subseteq I$ closed in $I$ and such that $[f_i,f_j]|_{X_I} = \id$
for all $i, j$.

This proposition will be proved by induction,
the key step being the following lemma.

\begin{lemma}
\label{lemma:step}
Let $f_1, \ldots, f_n$ be functions of class $C^2$,
$|f_i - \id|_{C^1} < \epsilon$
(where $\epsilon > 0$, $\epsilon < 1/100$).
Assume there exists a set $X \ne \emptyset$ closed in $(-1,1)$
and invariant under the functions $f_i$
such that $[f_i,[f_1,f_2]]|_X = \id$ for all $i$.
Then there exists $Y \subseteq X$, $Y \ne \emptyset$,
$Y$ also closed in $(-1,1)$ and invariant under the functions $f_i$
such that $[f_1,f_2]|_Y = \id$.
\end{lemma}

{\nobf Proof of Lemma \ref{lemma:step}:}
Let $Y$ be the set of fixed points of $[f_1,f_2]$:
the set $Y$ is clearly closed and invariant under the functions $f_i$;
it suffices to prove that $Y \ne \emptyset$.

If $X \cap (-1/2,1/2) = \emptyset$, let $x_0 \in X$
be the element of least absolute value;
assume without loss of generality that $x_0 > 0$.
We claim that $f_i(x_0) = x_0$ for all $i$:
indeed, $f_i(x_0) < x_0$ would be another element of $X$ of smaller
absolute value whence $f_i(x_0) \ge x_0$.
Similarly, $f_i^{-1}(x_0) \ge x_0$ and since $f'_i \ge 0$ the claim follows.
We can now take $Y = \{x_0\}$ and this proves the lemma in this case;
we assume from now on $X \cap (-1/2,1/2) \ne \emptyset$.

Let $x_0 \in X \cap (-1/2, 1/2)$ be an arbitrary point.
If $f_1(x_0) = f_2(x_0) = x_0$ we are done.
Otherwise we may assume without loss of generality that
$f_1^{-1}(x_0) < x_0 \le f_2(x_0) \le f_1(x_0)$, $x_0 < f_1(x_0)$.
We prove that there exists a fixed point
of the commutator $[f_1,f_2]$ in the interval $[f_1^{-1}(x_0),f_1(x_0)]$:
this will prove the lemma.

Set $\delta = f_1(x_0) - x_0 > 0$, $I_k = (f_1^{-k}(x_0), f_1^{k}(x_0))$.
Clearly $I_1 \subset I_2 \subset \cdots \subset I_{10} \subset (-3/4,3/4)$
and both $f_1(I_k)$ and $f_1^{-1}(I_k)$ are contained in $I_{k+1}$
for $k > 0$, $k < 10$.
Notice that $\delta (1-\epsilon)^{|k| + 1} <
f_1^{k+1}(x_0) - f_1^k(x_0) < \delta (1+\epsilon)^{|k| + 1}$
for $-10 \le k < 10$;
in particular, $f_1^k(x_0) - x_0 < k(1+\epsilon)^k \delta$
for $k > 0$, $k < 10$.

We claim that $f_2(I_k)$ and $f_2^{-1}(I_k)$ are both contained
in $I_{k+2}$ for all $k > 0$, $k < 9$.
From the mean value theorem,
$f_2f_1^k(x_0) - f_2(x_0) < (1+\epsilon)(f_1^k(x_0) - x_0)$
whence, adding $f_2(x_0) \le f_1(x_0)$,
$f_2f_1^k(x_0) < f_1^k(x_0) + (1 + \epsilon k (1 + \epsilon)^k)\delta$. 
On the other hand,
$f_1^{k+2}(x_0) > f_1^k(x_0) + 2(1 - \epsilon)^{k+1} \delta$
and it follows that $f_2f_1^k(x_0) < f_1^{k+2}(x_0)$
provided $2(1 - \epsilon)^{k+1} > 1 + \epsilon k (1+\epsilon)^k$,
which indeed happens for $\epsilon < 1/100$ and $k \le 8$.
Similarly, $f_2^{-1}f_1^k(x_0) < f_1^{k+2}(x_0)$,
$f_2f_1^{-k}(x_0) > f_1^{-k-2}(x_0)$ and
$f_2^{-1}f_1^{-k}(x_0) > f_1^{-k-2}(x_0)$,
proving our claim.

Let $g$ be the commutator $[f_1,f_2]$.
From the claim above it follows easily that 
$g(I_k)$ and $g^{-1}(I_k)$ are both contained
in $I_{k+4}$ for $k > 0$, $k < 6$.
For example,
$f_2f_1^k(x_0) < f_1^{k+2}(x_0)$ implies
$f_1f_2f_1^k(x_0) < f_1^{k+3}(x_0)$ and
$f_2^{-1}f_1f_2f_1^k(x_0) < f_1^{k+5}(x_0)$ whence, finally,
$[f_1,f_2] f_1^k(x_0) = f_1^{-1}f_2^{-1}f_1f_2f_1^k(x_0) < f_1^{k+4}(x_0)$.


Let $\psi: I_{10} \to \RR$
be a function with positive derivative of class $C^2$
satisfying $\psi(x_0) = 0$, $\psi(f_1(x)) = 1 + \psi(x)$:
such a function can be contructed by taking
a diffeomorphism from $[x_0,f_1(x_0)]$ to $[0,1]$
with compatible behavior at boundary points.
Define $\tilde f_i: (-8,8) \to \RR$ and $\tilde g: (-6,6) \to \RR$
by $\tilde f_i = \psi f_i \psi^{-1}$ and 
$\tilde g = \psi g \psi^{-1}$.
Notice that $\tilde f_1(x) = x+1$.
Also, set $\tilde X = \psi(X \cap I_{10})$:
if $x \in \tilde X, |x| < 6$
then the points $\tilde f_i(x)$, $\tilde f_i^{-1}(x)$,
$\tilde g(x)$ and $\tilde g^{-1}(x)$ are all in $\tilde X$.
Let $\check g: \RR \to \RR$ be a function of degree $1$
(i.e., $\check g(x+1) = \check g(x) + 1$ for all $x$)
with $\check g|_{\tilde X} = \tilde g|_{\tilde X}$.
In order to prove the existence of such a function $\check g$,
we consider two cases.
If $[0,1] \subseteq \tilde X$ then $\tilde g(x+1) = \tilde g(x) + 1$
for all $x \in (-5,5)$ and $\check g$ is obtained by extending $\tilde g$.
Otherwise, let $[x_1,x_2] \subset (0,1)$,
$[x_1,x_2] \cap \tilde X = \{x_1,x_2\}$
and define $\check g$ to be an arbitrary function of degree $1$
coinciding with $\tilde g$ in the interval $[x_2 - 1, x_1]$.
Let $c$ be the translation number of $\check g$:
we already proved that $|c| < 4$; we claim that $c = 0$.
The claim implies that the sequence
$0, \check g(0), \ldots, \check g^k(0), \ldots$
converges to a fixed point $\psi(x_1)$ of $\check g$
in $\tilde X \cap [-1,1]$ and $x_1$ is the desired fixed point of $g$,
proving the lemma.
In order to prove the claim,
it is convenient to consider, by contradiction,
the cases $c$ irrational and $c$ rational, $c \ne 0$.

{\nobf Case 1: $c$ irrational.}

By Denjoy theorem there exists a homeomorphism
$\phi: \RR \to \RR$ with $\phi(n) = n$ for $n \in \ZZ$ and 
\[\hat f_1(x) = \phi \tilde f_1 \phi^{-1}(x) = x + 1, \quad
\hat g(x) = \phi \check g \phi^{-1}(x) = x + c. \]
Define $\hat f_2: (-7,7) \to \RR$ by
\[ \hat f_2(x) =  \phi \tilde f_2 \phi^{-1}(x). \]
We have $\hat f_2 \hat f_1 \hat g = \hat f_1 \hat f_2$
and $\hat f_2 \hat g = \hat g \hat f_2$ on
$\hat X = \phi(\tilde X)$ which become
\[ \hat f_2(x + 1 + c) = \hat f_2(x) + 1, \qquad
\hat f_2(x + c) = \hat f_2(x) + c, \qquad x \in \hat X. \]
Also, since $f_2(I_k) \subset I_{k+2}$ we have
$\hat f_2((-k,k)) \subset (-k-2,k+2)$
for $k \le 7$.

Consider the set $A$ of points $(x,y) \in \ZZ^2$
with $|x(1+c) + yc| < 5$: this set is connected
in the sense that points of $A$ can be joined by a path
with vertices in $A$ and edges of size 1.
Let $(x_k, y_k)$, $k = 0, \ldots, N$ be such a path of points of $A$
with $(x_0, y_0) = (0,0)$ and $|x_N + y_N c| > 9$.
From $|x_k(1+c) + y_k c| < 5$ we have $|\hat f_2(x_k(1+c) + y_k c)| < 7$
for all $k \le N$.
If $(x_{k+1},y_{k+1}) = (\pm 1,0) + (x_k, y_k)$
then $\hat f_2(x_{k+1}(1+c) + y_{k+1}c) = \hat f_2(x_k (1+c) + y_k c) \pm 1$.
Also, if $(x_{k+1},y_{k+1}) = (0,\pm 1) + (x_k, y_k)$
then $\hat f_2(x_{k+1}(1+c) + y_{k+1}c) = \hat f_2(x_k (1+c) + y_k c) \pm c$.
Thus $\hat f_2(x_k(1+c) + y_k c) = \hat f_2(0) + x_k + y_k c$ for all $k$.
In particular $|x_N + y_N c| = |\hat f_2(x_N(1+c) + y_N c) - \hat f_2(0)| < 9$,
a contradiction.

\smallskip

{\nobf Case 2: $c$ rational, $c \ne 0$.}

Let $c = p/q$, $p, q \in \ZZ$, $q > 0$.
Assume without loss of generality that $p > 0$.
Let $h = \tilde f_1^{-p} \check g^q$:
we know that $h$ has a fixed point $x_1 \in [0,1]$.
Take $x_2 \in \tilde X \cap [x_1, x_1 + 1]$:
the sequence $x_2, h(x_2), h^2(x_2), \ldots$
is contained in the compact set $\tilde X \cap [x_1, x_1 + 1]$
and therefore any accumulation point $x_3$ of this sequence
is a fixed point of $h$ in $\tilde X$.

Define $z_0 = x_3$ and
\[ z_{i+1} = \begin{cases}
\tilde f_1^{-1}(z_i),&z_i \ge x_3 + 1,\\
\check g(z_i),&\textrm{otherwise}.
\end{cases} \]
We have $z_{p+q} = z_0$ and in this sequence we take
$p$ times the first case and $q$ times the second.
Set $w_i = \tilde f_2^{-1}(z_i)$:
we have therefore
\[ w_{i+1} = \begin{cases}
\check g^{-1} \tilde f_1^{-1}(w_i),&w_i \ge f_2^{-1}(x_3 + 1),\\
\check g(w_i),&\textrm{otherwise}.
\end{cases} \]
Thus $\tilde f_1^{-p} \check g^{q-p}$ has fixed point $w_0$
and the translation number of $\check g$ is $p/(q-p)$,
a contradiction. \qed



\bigskip

{\nobf Proof of Proposition \ref{prop:fixedpoint}:}
We proceed by induction on $m$.
Assume our pseudogroup of functions to be nilpotent of order $m$.
Apply lemma \ref{lemma:step} to the family
of commutators of order at most $m-1$, $X_0 = (-1,1)$,
with the new $f_1$ being an arbitrary commutator of order at most $m-1$
and the new $f_2$ being a commutator or order $m-1$.
We thus obtain a closed invariant subset $X_1$ of $X_0$
where a given commutator of order $m$ equals the identity.
Repeating this process we obtain a closed invariant subsets
$X_2 \supseteq X_3 \supseteq \cdots \supseteq X_k$
such that all commutators of order $m$ equal the identity in $X_k$.
Now apply the induction hypothesis to obtain $X \subseteq X_k$.
\qed

If in the constructions performed in the proof of the lemma we take $x_0$
to be a fixed point of $g$ then $f_1f_2(x_0) = f_2f_1(x_0)$ and
$x_0 \le f_2(x_0) \le f_1(x_0)$ implies
$f_1(x_0) \le f_2f_1(x_0) = f_1f_2(x_0) \le f_1^2(x_0)$
and $f_2(I_k) \subseteq I_{k+1}$. Also, $g(I_k) = I_k$.

\section{Translation number}

We modify the definition of translation number to define the translation
number of $f_2$ relative to $f_1$, where both functions
belong to a near-identity nilpotent pseudogroup.
Let $x_0$ be a point of $X$,
a nontrivial invariant set as discussed in the previous section.
Assume that $f_1^{-1}(x_0) < x_0 \le f_2(x_0) \le f_1(x_0)$.
Define a sequence of points starting at $a_0 = x_0$ by 
$a_{n+1} = f_1^{-k(n)}f_2(a_n)$ where $k(n)$ is the only integer
for which $x_0 \le f_1^{-k(n)}f_2(a_n) < f_1(x_0)$.
Since $f_1(x_0) \le f_1f_2(x_0) = f_2f_1(x_0) \le f_1^2(x_0)$
and, by contruction, $x_0 \le a_n < f_1(x_0)$
we have $x_0 \le f_2(x_0) \le f_2(a_n) < f_2f_1(x_0) \le f_1^2(x_0)$
and therefore $k(n)$ is $0$ or $1$.
Also, $k(n) = 1$ if and only if $x_0 \le a_{n+1} < f_2(x_0)$.
Let $F^{x_0}_0 = \id$ and
$F^{x_0}_{n+1} = f_1^{-k(n)} \circ f_2 \circ F^{x_0}_n$:
as usual, we must show that this composition makes sense in a large domain.
As in the proof of proposition \ref{prop:fixedpoint},
let $I_n = (f_1^{-n}(x_0), f_1^n(x_0))$;
notice that $I_n$ is well defined at least for $n \le 10$
and that $f_2(I_n) \subseteq I_{n+1}$.

\begin{lemma}
\label{lemma:goodFn}
For any positive integer $n$, $F^{x_0}_n$ is well defined in $I_8$
and $F^{x_0}_n(I_8) \subseteq I_9$. Also,
\[ f_1^{-8}(x_0) \le F^{x_0}_n f_1^{-8}(x_0) =
f_1^{-8} F^{x_0}_n(x_0) < f_1^{-7}(x_0),\]
\[ f_1^{8}(x_0) \le F^{x_0}_n f_1^{8}(x_1) =
f_1^{8} F^{x_0}_n(x_1) < f_1^{9}(x_0).\]
Moreover, $X$ is invariant under $F^{x_0}_n$.
\end{lemma}

\proof
We prove the inequalities in the statement by induction on $n$,
the case $n = 1$ being easy.
By definition,
$F_{n+1}f_1^{\pm 8}(x_0) = f_1^{-k(n)}f_2F_nf_1^{\pm 8}(x_0)$:
this already shows that these two expressions make sense
since by induction hypothesis both $F_nf_1^{\pm 8}(x_0)$
make sense and are in $I_9$.
By the induction hypothesis,
$F_{n+1}f_1^{\pm 8}(x_0) = f_1^{-k(n)}f_2f_1^{\pm 8}F_n(x_0)$.
Since $F_n(x_0) \in X$,
$F_{n+1}f_1^{\pm 8}(x_0) = f_1^{\pm 8}f_1^{-k(n)}f_2F_n(x_0) =
f_1^{\pm 8}F_{n+1}(x_0)$.
The other claims are now easy.
\qed

Set $p(0) = 0$, $p(n+1) = p(n) + k(n)$:
we define the translation number to be the limit
\[ \rot(f_2,f_1,x_0) = \lim_{n \to \infty}{\frac{p(n)}{n}}; \]
we still have to prove that this limit exists.
If $f_2^{-1}(x_0) < x_0 \le f_1(x_0) \le f_2(x_0)$
we can make a similar construction reverting the roles
of $f_1$ and $f_2$ and define
\[ \rot(f_2,f_1,x_0) = 1/(\rot(f_1,f_2,x_0)); \]
if $f_1^{-1}(x_0) \le f_2(x_0) \le x_0 < f_1(x_0)$
we define
\[ \rot(f_2,f_1,x_0) = -\rot(f_2^{-1},f_1,x_0); \]
the other cases are similar.

Let $Z$ be the closed set of common fixed points of all functions $f_i$.
We show that the translation number is well defined in each connected
component of the complement of $Z$.

\begin{prop}
\label{prop:translationnumber}
If $x_0, x_1 \in X$ are in the same connected component
of the complement of $Z$ then
$\rot(f_2,f_1,x_0) = \rot(f_2,f_1,x_1)$.
\end{prop}

\proof
Let $g = [f_1,f_2]$. We first prove that the limit exists.
Assume $f_1^{-1}(x_0) < x_0 \le f_2(x_0) \le f_1(x_0)$
and let $\psi$ be a conjugation between $f_1$ and $x \mapsto x+1$
so that $\tilde f_1(x) = \psi f_1 \psi^{-1}(x) = x+1$,
$\tilde f_2 = \psi f_2 \psi^{-1}$ and $\tilde g = \psi g \psi^{-1}$;
$\psi: I_{10} \to \RR$ is a function with positive derivative of class $C^2$
and $\psi(x_0) = 0$ as constructed in the proof
of Lemma \ref{lemma:step}.
Notice that $\tilde f_2(x + 1) = \tilde f_2(x) + 1$
if and only if $\tilde g(x) = x$;
in particular $\tilde f_2(1) = \tilde f_2(0) + 1$.
The definition of translation number applies to the restriction
to $[0,1]$ of $\tilde f_2$;
equivalently, let $\check f_2$ be the only function of degree $1$
coinciding with $\tilde f_2$ in the interval $[0,1]$:
\[ \check f_2(x) = \tilde f_2(x - \lfloor x \rfloor) + \lfloor x \rfloor. \]
The points $a_n$ constructed as above from $f_1$ and $f_2$
are all in the interval $[x_0, f_1(x_0))$
and $\psi(a_n)$ is always in the interval $[0,1]$.
The construction of $k(n)$, $F_n$ and $p(n)$
only considers values of $f_2$ in the interval $[x_0, f_1(x_0))$,
or, equivalently, values of $\tilde f_2$ in the interval $[0,1)$.
It makes therefore no difference
whether we take $\tilde f_2$ or $\check f_2$
and $\rot(f_2,f_1,x_0)$ is the usual translation number of $\check f_2$.

Let $x_1$ be another fixed point of $g$, $x_0 \le x_1 \le f_1(x_0)$.
The translation number of $\check f_2$ is the same if computed in the interval
$[0,1]$ or in the interval $[\psi(x_1), \psi(x_1) + 1]$.
Furthermore, the functions $\tilde f_2$ and $\check f_2$ coincide in the orbit
of $x_1$ since these are all fixed points of $\tilde g$
and the construction of the translation number coincides for these two
functions. Thus $\rot(f_2,f_1,x_0) = \rot(f_2,f_1,x_1)$.

Let $x_0 < x_\ast$ be two fixed points of $g$ in the same connected
component of the complement of $Z$. Let $\epsilon > 0$ be the infimum
over the compact interval $[x_0, x_\ast]$ of the positive continuous
function $\max\{ |f_1(x) - x|, |f_2(x) - x| \}$ and take a sequence
$y_0 = x_0, y_1, \ldots, y_N = x_\ast$ with $0 < y_{i+1} - y_i < \epsilon/4$.
Let $x_i$ be the fixed point of $g$ which is closest to $y_i$
so that $x_0 \le x_1 \le \cdots \le x_N = x_\ast$ are fixed points of $g$.
We claim that
$x_{i+1} \le \max \{ f_1(x_i), f_2(x_i), f_1^{-1}(x_i), f_2^{-1}(x_i) \}$.
Assume without loss of gererality that the largest among these
four numbers is $f_1(x_i)$: the interval $((x_i + f_1(x_i))/2, f_1(x_i))$
has size at least $\epsilon/2$ and therefore there is some $y_j$ in it.
The point $f_1(x_i)$ is a fixed point of $g$ thus the distance between
$x_j$ and $y_j$ is no larger than that between $f_1(x_i)$ and $y_j$:
it follows that $x_i \le x_{i+1} \le x_j \le f_1(x_i)$, as claimed.
The previous paragraph can now be used to show that
$\rot(f_2,f_1,x_i) = \rot(f_2,f_1,x_{i+1})$ and we have
$\rot(f_2,f_1,x_0) = \rot(f_2,f_1,x_N)$.
\qed

\section{The irrational case}

We already saw that commutators in a near-identity nilpotent
pseudogroup of functions have many fixed points;
we now show that in many cases
all commutators equal the identity so that the original functions commute.

\begin{prop}
\label{prop:irrational}
Let $f_1, f_2$ be functions in a near-identity nilpotent pseudogroup
and let $X$ be a closed invariant subset where all the $f_i$ commute.
Let $x_0 \in X$, $x_0$ not a fixed point of $f_1$.
If $\rot(f_2,f_1,x_0)$ is irrational then
$x_0$ is an interior point of $X$.
\end{prop}

\proof
Assume without loss of generality that
$f_1^{-1}(x_0) < x_0 \le f_2(x_0) \le f_1(x_0)$.
Construct $\tilde f_1$, $\tilde f_2$, $\tilde g$ and $\tilde X$
as usual where $g = [f_1,f_2]$ so that $\tilde f_1(x) = x+1$;
notice that $\tilde g(0) = 0$.
Assume by contradiction that $(x_1, x_2) \subseteq (0, 1)$ is an open interval
such that $[x_1, x_2] \cap \tilde X = \{x_1,x_2\}$.
We shall construct a counter-example to Denjoy's theorem,
thus obtaining a contradiction.

Let $\check f_2: \RR \to \RR$ be an increasing function of degree 1
and class $C^2$ with $\check f_2(x) = \tilde f_2(x)$ for $x_2 - 1 < x < x_1$.
The function $\check f_2$ is defined arbitrarily in the interval $(x_1, x_2)$
with the only restrictions that it must be of class $C^2$, increasing,
satisfy $\check f_2(x) = \tilde f_2(x)$ for $x$ near $x_1$
and $\check f_2(x) = 1 + \tilde f_2(x-1)$ for $x$ near $x_2$;
this is clearly possible.
Now $\check f_2$ is a function of degree 1 and class $C^2$
and irrational translation number and $\tilde X + \ZZ$
is a nontrivial invariant closed set,
contradicting Denjoy's theorem.
\qed

We may bring together our conclusions as a proposition.

\begin{prop}
\label{prop:3cases}
Let $f_1, \ldots, f_n$ be a near-identity nilpotent pseudogroup of functions
and let $x_0 \in (-1/2,1/2)$. Then one of the following three
situations holds:
\begin{enumerate}
\item{$x_0$ is a common fixed point of the functions $f_i$.}
\item{There exists real constants $a_i$,
an open interval $I \subseteq (-1,1)$ 
containing $x_0, f_i(x_0), f_i^{-1}(x_0)$,
and a homemorphism $\phi: J \to I$, $J \subseteq \RR$ with 
$f_i(\phi(t)) = \phi(t + a_i)$ whenever $t, t+a_i \in J$.}
\item{There exist integer constants $a_i$,
a finite set $\{y_{-N}, y_{-N+1}, \ldots, y_N\}$ with $y_i < y_{i+1}$
and $y_0 \le x_0 < y_1$, $N > |a_i|$ with
$f_i(y_k) = y_{k + a_i}$.}
\end{enumerate}
\end{prop}

\proof
If $x_0$ is not a common fixed point of the functions $f_i$
then we may apply the results of the previous sections.
If at least one translation number is irrational then
we apply proposition \ref{prop:irrational} and we are in the second case.
Otherwise the functions from $\Ss^1$ to itself constructed above
all have rational translation numbers and therefore admit periodic points
and we are in the third case.
\qed

\section{Near-identity nilpotent pseudogroups \\ are metabelian}

We first state Koppel's lemma (\cite{K}), an important result
also in the works of Plante, Thurston, Farb and Franks.

\begin{lemma}
\label{lemma:koppel}
Let $f, g: \RR \to \RR$ be increasing diffeomorphisms,
$f$ of order $C^2$ and $g$ of order $C^1$ with $[f,g] = \id$.
Assume that there exists a nondegenerate bounded open interval $I = (a,b)$
for which $f(b) = g(b) = b$, $f(a) > g(a) = a$
and $x \in (a,b)$ implies $f(x) > x$. Then $g(x) = x$ for all $x \in I$.
\end{lemma}

Let $f_1, f_2, \ldots, f_n$ be a near-identity nilpotent pseudogroup
with rational
translation numbers as in item 3 of proposition \ref{prop:3cases}:
let $y_0, y_1, a_1, a_2, \ldots, a_n$ be as in the proposition.
We define a second near-identity nilpotent pseudogroup
$\tilde f_1, \ldots, \tilde f_N$
of functions satisfying $\tilde f_i(y_0) = y_0$, $\tilde f_i(y_1) = y_1$,
$N = \binom{n}{2} + n - 1$.

Assume without loss of generality that $a_n > 0$.
The construction of $\tilde f_i$, $i < n$, is similar to that
of $F^{x_0}_j$ in section 4: define $F_{i,0} = \id$ and
\[ F_{i,j+1} = f_n^{-\lfloor k/a_n \rfloor} \circ f_i \circ F_{i,j}
\quad \textrm{where} \quad f_i(F_{i,j}(y_0)) = y_k. \]
In this way $F_{i,j}(y_0) \ge y_0$ for all $i, j$.
We define $\tilde f_i = F_{i,a_n}$:
we have $\tilde f_i(y_0) = y_0$ and $\tilde f_i(y_1) = y_1$.
The remaining $\tilde f_i$ are the commutators $[f_j,f_k]$, $1 \le j < k \le n$.
We claim that the functions $\tilde f_i$ commute:
this will establish our claim that the original pseudogroup is metabelian.

Apply proposition \ref{prop:3cases} to the pseudogroup $\tilde f_i$
at some point in the interval $(y_0, y_1)$:
if we have cases 1 or 2 we are done.
We assume therefore that we have case 3
in a maximal interval $(a,b) \subseteq (y_0, y_1)$.
Let $\tilde y_0, \tilde y_1, \tilde a_1, \ldots, \tilde a_N$
be as in proposition \ref{prop:3cases}.
Assume without loss of generality that $\tilde a_k > 0$
so that $\tilde f_k(a) = a$, $\tilde f_k(b) = b$
and $\tilde f_k(x) > x$ for all $x \in (a,b)$.
Assuming that the functions $\tilde f_i$ do not commute,
let $g$ be a commutator of highest order which is different
from the identity. By construction $[\tilde f_k, g] = \id$.
If $g$ is not one of the $\tilde f_i$ then $g(\tilde y_0) = \tilde y_0$
so that $g$ has a fixed point in $(a,b)$. Lemma \ref{lemma:koppel} now
implies $g = \id$ which is a contradiction unless $g = \tilde f_i$
and the functions $\tilde f_i$ commute, as required.

\section{Proof of theorems \ref{theo:classify} and \ref{theo:foliation}}

We proved in the previous section
that any near-identity pseudogroup is metabelian.

If at least one common fixed point exists, we are in case 1.
Let $I_1 = (a,c)$ be a maximal open interval containing no common
fixed points. Assume without loss of generality that $|a| < 1$
and apply proposition \ref{prop:fixedpoint} to the interval $(a,c)$ to obtain
$b \in (a,c)$, a fixed point for all commutators in the pseudogroup.
Let $f = f_j$ or $f = f_j^{-1}$ be such that $f(b)$ is maximal;
from the results of section 3, $f(x) > x$ for all $x \in (a,c)$.
Take $g$ to be an arbitrary commutator and apply Koppel's lemma
(\ref{lemma:koppel}) to obtain $g(x) = x$ for all $x \in (a,c)$.
This implies that the pseudogroup is abelian.
For each $i$, let $a_i =  \rot(f_i,f,b)$.
If at least one of the $a_i$ is irrational, the existence of the
homeomorphism $\phi$ follows from Denjoy's theorem.
Otherwise, we may write $a_i = p_i/q$, where $p_i$ and $q$ are integers.
Assume without loss of generality that $p_k = 1$ so that, for each $i$,
$\rot(f_i f_k^{-p_i}, f, b) = 0$:
this implies the existence of $\tilde b_i \in [b,f(b)]$ with
$f_i f_k^{-p_i}(\tilde b_i) = \tilde b_i$.
Koppel's lemma now implies that $f_i = f_k^{p_i}$
and the pseudogroup consists of integer powers of $f_k$,
implying the existence of the homeomorphism $\phi$.
This concludes the discussion of case 1.

If no common fixed point exists,
take $x_0 \in X$ and let $f = f_j$ or $f = f_j^{-1}$ be such that
$f(x_0)$ is maximal and set $a_i = \rot(f_i,f,x_0)$.
If at least one of the $a_i$ is irrational
then from section 4 the pseudogroup is abelian and from Denjoy's theorem
there exists a homeomorphism $\phi$ as stated, concluding case 2.
Otherwise we are in case 3 of proposition \ref{prop:3cases},
concluding theorem \ref{theo:classify}.

For theorem \ref{theo:foliation}, take the holonomy pseudogroup and apply
theorem \ref{theo:classify}: each case in one theorem corresponds to the
same case in the other. In cases 1 and 2, the homomorphism $\alpha: G \to \RR$
takes each generator $\gamma_i$ of $\pi_1(G/H)$ to $a_i$ as constructed
in theorem \ref{theo:classify} where $f_i = f_{\gamma_i}$.
The homeomorphism $\Phi$ is constructed from $\phi$.
More precisely,
for $(gH,x) \in G/H$ let $\gamma: [0,1] \to G$ be a path with
$\gamma(0) \in H$ and $\gamma(1) = g$.
Lift $\gamma$ to $\Ff_\alpha$ to obtain a path
$\tilde\gamma: [0,1] \to G/H \times \RR$ tangent to $\Ff_\alpha$
with $\tilde\gamma(1) = (gH,x)$: let $\tilde\gamma(0) = (H,x_0)$.
Now lift $\gamma$ to $\Ff_1$ to obtain $\hat\gamma: [0,1] \to G/H \times (-1,1)$
with $\hat\gamma(0) = (H,\phi(x_0))$ and define $\Phi(gH,x) = \hat\gamma(1)$.
The properties of $\phi$ imply that $\Phi(gH,x)$ is well defined, i.e.,
does not depend on the choice of $\gamma$.
This concludes the proof of theorem \ref{theo:foliation}.

\vfil\eject

{\obeylines
\parskip=0pt
\parindent=0pt

Tania M. Begazo
Universidade Federal Fluminense
Instituto de Matem\'atica
Niterói, RJ 24020-140, Brazil
tania@mat.uff.br

\smallskip

Nicolau C. Saldanha
Depto. de Matem{\'a}tica, PUC-Rio
R. Mq. de S. Vicente 225
Rio de Janeiro, RJ 22453-900, Brazil
nicolau@mat.puc-rio.br
http://www.mat.puc-rio.br/$\sim$nicolau
}

\end{document}